\documentclass{amsart}
\usepackage{amssymb}
\usepackage{pdflscape}
\usepackage{amscd}
\usepackage{fancyhdr}
\usepackage{graphicx}

\usepackage[usenames,dvipsnames]{xcolor} 
\usepackage{tikz, tikz-3dplot}
\usepackage{tkz-graph}
\usepackage{tikz-cd}
\usetikzlibrary{positioning, calc, arrows, decorations.markings, decorations.pathreplacing}
\tikzset{
  symbol/.style={
    draw=none,
    every to/.append style={
      edge node={node [sloped, allow upside down, auto=false]{$#1$}}}
  }
}
\tikzcdset{scale cd/.style={every label/.append style={scale=#1},
    cells={nodes={scale=#1}}}}

\newtheorem{theorem}{Theorem}[section]
\newtheorem{corollary}[theorem]{Corollary}
\newtheorem{lemma}[theorem]{Lemma}
\newtheorem{proposition}[theorem]{Proposition}

\newtheorem{definition}[theorem]{Definition}

\begin{document}
\title{Expander representations of quivers}
\author{Markus Reineke}
\begin{abstract} We propose a definition of expander representations of quivers, generalizing dimension (or linear algebra) expanders, as a qualitative refinement of slope stability. We prove existence of uniform expander representations for any wild quiver over an algebraically closed base field, using the concept of general subrepresentations and spectral properties of Cartan matrices.
\end{abstract}
\maketitle
\parindent0pt

\section{Introduction}

Dimension expanders were first introduced in \cite{LZ}, as a linear algebra analogue of the highly successful notion of expander graphs; see \cite{Lu} for a survey of the many fascinating applications of expanders. The relation between different graph-theoretic and linear algebra notions of expansion is explained in detail in \cite{WW}.\\[1ex]
In \cite{REx}, dimension expanders are interpreted in terms of stable representations of generalized Kronecker quivers, and results on general subrepresentation types are used to give a new proof of existence of expanders, as well as to provide optimal bounds for the expansion coefficients.\\[1ex]
Quiver representations, initially a language for encoding classification problems of linear algebra, allow to enhance such classification problems with techniques of homological algebra, by forming a hereditary abelian category, and of algebraic geometry, by forming orbits for algebraic group actions. It is exactly this interplay which yields criteria for (non-)existence of subrepresentations \cite{Scho}, leading to the new existence proof for dimension expanders in \cite{REx}. Surprisingly, these quiver-theoretic techniques are strong enough to even yield optimal expansion coefficients, using special properties of the generalized Kronecker quivers.\\[1ex]
With this in mind, it is thus desirable to explore the utility of quiver representation theory for expander phenomena further, and in particular to ask for the quiver-theoretic nature of explicit expansion coefficients. Abstracting from \cite{REx}, one can ask for a generalization of the concept of dimension expanders using a quantitative version of slope stability; such a generalization is briefly discussed in \cite[Definition 5.2]{REx}.\\[1ex]
In the present work, we elaborate on this and propose a definition of (uniform) expanders for representations of quivers (Definitions \ref{defexp}, \ref{defuniexp}). Our main result, Theorem \ref{main}, states that every wild quiver admits a uniform expander for some slope function, and this already characterizes wildness of the quiver. Using general subrepresentation types as in \cite{REx}, we reduce the existence proof to spectral properties of the Cartan matrix of the quiver. More precisely, a uniform expansion coefficient ultimately can be constructed using the spectral gap of the Cartan matrix, that is, the difference between its two smallest eigenvalues -- in full agreement with the role of spectral gaps in expander theory.\\[1ex]
After fixing some quiver terminology, we introduce our concept of (uniform) expander in Section \ref{s2}. We give a numerical criterion for their existence in Section \ref{s3}, using results of \cite{Scho} on general subrepresentations. We show in Section \ref{sectionkronecker} that the present concept is compatible with the methods of \cite{REx} in the case of generalized Kronecker quivers, providing explicit expansion coefficients in this case. In the central Section \ref{wild}, we prove existence of uniform expanders for arbitrary wild quivers using spectral properties of the Cartan matrix. 
That uniform expanders require the quiver to be wild is proven in Section \ref{nonpp}, where we, more generally, prove that uniform expanders cannot consist of preprojective/preinjective representations only. Finally we prove our main result in Section \ref{concluding}, and conclude with some remarks.\\[2ex]
{\bf Acknowledgments:} The author is indebted to Urban Jezernik for generously sharing a crucial argument used in the proof of Lemma \ref{boundbydd}, which allowed to simplify the proof in comparison to an earlier version of this work. The author would like to thank Sebastian Eckert, Lutz Hille and Alastair King for several helpful discussions, and an anonymous referee for helpful suggestions.

\section{Expander representations and uniform expansion}\label{s2}

For all general notions about quivers and their representations, we refer to \cite{Schiffler}. 
Let $Q$ be a finite quiver, given by a finite set $Q_0$ of vertices and finitely many arrows written $\alpha:i\rightarrow j$. We assume $Q$ to be acyclic, that is, $Q$ has no oriented cycles. We consider the real vector space $\mathbb{R}Q_0$ with coordinate basis vectors ${\bf i}$ for $i\in Q_0$, and write vectors ${\bf d}\in \mathbb{R}Q_0$ as
$${\bf d}=\sum_{i\in Q_0}d_i{\bf i}.$$ 
We define the Euler form of $Q$ as the bilinear form
$$\langle{\bf d},{\bf e}\rangle=\sum_{i\in Q_0}d_ie_i-\sum_{\alpha:i\rightarrow j}d_ie_j.$$
We denote by $(\_,\_)$ its symmetrization, and by $\{\_,\_\}$ its antisymmetrization. With respect to the coordinate basis, $(\_,\_)$ is represented by a matrix $C$ which is a symmetric generalized Cartan matrix.\\[1ex]
We fix an algebraically closed base field $F$ and consider finite-dimensional $F$-representations $V$ of $Q$. Thus, $V$ is given by finite-dimensional $F$-vector spaces $V_i$ for $i\in Q_0$ and $F$-linear maps $V_\alpha:V_i\rightarrow V_j$ for all arrows $\alpha:i\rightarrow j$. Morphisms between representations $V$ and $W$ are tuples $(\varphi_i:V_i\rightarrow W_i)_{i\in Q_0}$ of linear maps intertwining with the structure maps, that is, 
$$W_\alpha\circ \varphi_i=\varphi_j\circ V_\alpha$$
for all arrows $\alpha:i\rightarrow j$ in $Q$. Defining composition of morphisms componentwise results in an abelian $F$-linear category ${\rm rep}_FQ$ of finite-dimensional $F$-representations of $Q$, which is equivalent to the category of finite dimensional left modules over the path algebra $FQ$ of $Q$. In particular, we can apply all standard concepts of the theory of modules over rings to quiver representations. The dimension vector of a representation $V$ is the vector $${\rm\bf dim}(V)=\sum_{i\in Q_0}(\dim_FV_i){\bf i}\in\mathbb{R}_{\geq 0}Q_0.$$

We next introduce stability and refer to \cite{HD, K, Ru} for generalities. Let $\Theta,\kappa\in(\mathbb{R}Q_0)^*$ be two linear functionals, such that $\kappa$ is positive in the sense that $$\kappa(\mathbb{R}_{\geq 0}Q_0\setminus\{0\})\subset\mathbb{R}_{>0}$$ (equivalently, $\kappa({\bf i})>0$ for all $i\in Q_0$). We define the associated slope function 
$$\mu({\bf d})=\Theta({\bf d})/\kappa({\bf d})$$
for ${\bf d}\in\mathbb{R}_{\geq 0}Q_0\setminus\{0\}$.
We extend the functional $\Theta$ to representations by
$$\Theta(V)=\Theta({\rm\bf dim}(V)),$$
and similarly for $\kappa$ and for the slope function $\mu$.\\[1ex]
We call a non-zero representation $V$ stable with respect to the slope function $\mu$, or shortly $\mu$-stable, if $\mu(U)<\mu(V)$ for all non-zero proper subrepresentations $U\subset V$.  Note that the class of stable representations does not change when $\Theta,\kappa$ are replaced by positive scalar multiples, or when $\Theta$ is shifted by a multiple of $\kappa$. Stable representations $V$ are Schurian, that is, they have trivial endomorphism ring ${\rm End}_Q(V)\simeq F$, by a Schur's lemma type argument.\\[1ex]
From now one, we will fix a slope function $\mu$. We can now proceed to our main new concept. 

\begin{definition}\label{defexp} For $0<\delta<1$ and $\epsilon>0$, a representation $V$ of $Q$ is called a $(\delta,\epsilon)$-expander (with respect to $\mu$) if for all non-zero subrepresentations $U\subset V$ such that $\kappa(U)\leq\delta\cdot\kappa(V)$, we have $\mu(U)<\mu(V)-\epsilon$.
\end{definition}

The precise relation to stability is easily derived:

\begin{lemma} The following are equivalent for a non-zero representation $V$:
\begin{enumerate}
\item $V$ is $\mu$-stable,
\item for all $0<\delta<1$, there exists $\epsilon=\epsilon(\delta)>0$ such that $V$ is a $(\delta,\epsilon)$-expander (with respect to $\mu$).
\end{enumerate}
In particular, an expander representation is Schurian.
\end{lemma}

\proof Assume that the second condition is fulfilled, and let $U\subset V$ be a non-zero proper subrepresentation. Then $\kappa(U)=\delta\cdot\kappa(V)$ for some $\delta\in]0,1[$, thus 
$$\mu(U)\leq\mu(V)-\epsilon(\delta)<\mu(V),$$
thus $V$ is $\mu$-stable. Conversely, assume that $V$ is $\mu$-stable, and fix $\delta\in]0,1[$. If $V$ is simple, there is nothing to prove. Otherwise, there are only finitely many dimension vectors ${\bf e}$ of subrepresentations of $V$ such that $\kappa({\bf e})=\delta\cdot\kappa(V)$, thus
$$\epsilon(\delta)=\min\{\mu(V)-\mu(U)\, |\, 0\not=U\subset V,\, \kappa(U)\leq\delta\cdot\kappa(V)\}$$
exists, and is positive by $\mu$-stability of $V$. Thus, the second condition is fulfilled, finishing the proof.\\[1ex]
More important than individual expander representations are families of expanders for which the bound $\epsilon=\epsilon(\delta)$ can be chosen uniformly for the whole family:

\begin{definition}\label{defuniexp} A family of representations $(V^{(k)})_{k\geq 1}$ of $Q$ is called unbounded if for all $d$, there exists an index $k(d)$ such that $\dim V^{(k(d))}\geq d$.

An unbounded family is called a uniform $(\delta,\epsilon)$-expander if every $V^{(k)}$ is a $(\delta,\epsilon)$-expander. It is called a uniform expander if for every $\delta\in]0,1[$, there exists an $\epsilon=\epsilon(\delta)>0$ such that it is a uniform $(\delta,\epsilon)$-expander.

An unbounded sequence of dimension vectors $({\bf d}^k)_k$ is said to support a uniform expander if there exists a uniform expander $(V^{(k)})_k$ such that ${\rm\bf dim}(V^{(k)})={\bf d}^k$ for all $k\geq 1$. Finally, we say that the quiver $Q$, together with the slope function $\mu$, exhibits uniform expansion, if it admits a uniform expander.
\end{definition}

Our main result is:

\begin{theorem}\label{main} Every wild quiver exhibits uniform expansion for some slope function. Conversely, any such quiver is wild.
\end{theorem}

This will be proved in Section \ref{concluding}.

\section{Numerical criterion}\label{s3}

We first recall methods and results from \cite{Scho}.\\[1ex]
For two dimension vectors ${\bf e}\leq{\bf d}$ of representations of $Q$, we write ${\bf e}\hookrightarrow{\bf d}$ if every representation of dimension vector ${\bf d}$ admits a subrepresentation of dimension vector ${\bf e}$. This condition can be characterized numerically, although recursively \cite{CB,Scho}:

\begin{theorem} We have ${\bf e}\hookrightarrow{\bf d}$ if and only if $\langle{\bf e}',{\bf d}-{\bf e}\rangle\geq 0$ for all ${\bf e}'\hookrightarrow{\bf e}$.

In particular, ${\bf e}\hookrightarrow{\bf d}$ always implies $\langle{\bf e},{\bf d}-{\bf e}\rangle\geq 0$.
\end{theorem}

For a given dimension vector ${\bf d}$, and fixed vector spaces $V_i$ of dimension $d_i$ for $i\in Q_0$, respectively, the representations $V$ of $Q$ on the spaces $V_i$ are naturally parametrized by the affine space $$R_{\bf d}(Q)=\bigoplus_{\alpha:i\rightarrow j}{\rm Hom}_F(V_i,V_j),$$
which we consider with its Zariski topology. For another dimension vector ${\bf e}\leq{\bf d}$, the subset $S_{\bf e}\subset R_{\bf d}(Q)$ of representations admitting a subrepresentation of dimension vector ${\bf e}$ is known to be Zariski-closed. By definition, we have ${\bf e}\hookrightarrow{\bf d}$ if and only if $S_{\bf e}=R_{\bf d}(Q)$.\\[1ex]
The following is the key to existence proofs for (uniform) expanders. The proof is non-constructive and requires $F$ to be algebraically closed, but it effectively reduces existence to a purely numerical/root-theoretic problem.

\begin{lemma}\label{lemmaexistence}  Given a dimension vector ${\bf d}$ and $\delta,\epsilon$ as above, there exists a $(\delta,\epsilon)$-expander of dimension vector ${\bf d}$ if and only if, for all $0\not={\bf e}\hookrightarrow{\bf d}$ such that $\kappa({\bf e})\leq\delta\cdot\kappa({\bf d})$, we have
$\mu({\bf e})\leq\mu({\bf d})-\epsilon$.
\end{lemma}

\proof Let $V$ be a $(\delta,\epsilon)$-expander of dimension vector ${\bf d}$, and let $0\not={\bf e}\hookrightarrow{\bf d}$ be a dimension vector such that $\kappa({\bf e})\leq\delta\cdot\kappa({\bf d})$. Then there exists a subrepresentation $U\subset V$ of dimension vector ${\bf e}$, for which $\kappa(U)\leq \delta\cdot\kappa(V)$, and thus $$\mu({\bf e})=\mu(U)\leq\mu(V)-\epsilon=\mu({\bf d})-\epsilon.$$
Conversely, assume that the numerical condition holds. Consider the (finite) set $D$ of all dimension vectors $0\not={\bf e}\leq{\bf d}$ such that $\kappa({\bf e})\leq\delta\cdot\kappa({\bf d})$, but $\mu({\bf e})>\mu({\bf d})-\epsilon$. For every ${\bf e}$ in $D$, the condition ${\bf e}\hookrightarrow{\bf d}$ does not hold by assumption, thus $S_{\bf e}\subset R_{\bf d}(Q)$ is a proper Zariski-closed subset. Thus, the finite union
$$\bigcup_{{\bf e}\in D}S_{\bf e}\subset R_{\bf d}(Q)$$
is still a proper closed subset, and we can choose a representation $V$ in its complement. If $U\subset V$ is a subrepresentation such that $\kappa(U)\leq\delta\cdot\kappa(V)$, and ${\bf e}={\rm\bf dim}(U)$, then ${\bf e}$ does not belong to $D$ by the choice of $V$. Thus
$$\mu(U)=\mu({\bf e})\leq\mu(V)-\epsilon,$$
proving that $V$ is a $(\delta,\epsilon)$-expander.\\[1ex]
This lemma motivates the following notation:

\begin{definition} Given a dimension vector ${\bf d}$ and $\delta\in]0,1[$, define
$$\epsilon_{\bf d}^{\rm opt}(\delta)=\min\{\mu({\bf d})-\mu({\bf e})\, |\, 0\not={\bf e}\hookrightarrow{\bf d},\, \kappa({\bf e})\leq\delta\cdot\kappa({\bf d})\}$$
and
$$\epsilon_{\bf d}^{\rm eff}(\delta)=\min\{\mu({\bf d})-\mu({\bf e})\, |\, 0\not={\bf e}\leq{\bf d},\, \langle{\bf e},{\bf d}-{\bf e}\rangle\geq 0,\, \kappa({\bf e})\leq\delta\cdot\kappa({\bf d})\}.$$
\end{definition}

Thus, by the previous lemma, $\epsilon_{\bf d}^{\rm opt}(\delta)$ is the largest $\epsilon$ for which a $(\delta,\epsilon)$-expander of dimension vector ${\bf d}$ exists. It is difficult to determine, though, due to the recursive nature of the criterion for deciding ${\bf e}\hookrightarrow{\bf d}$ (see, however, Section \ref{sectionkronecker} for the case of generalized Kronecker quivers). In contrast, since ${\bf e}\hookrightarrow{\bf d}$ implies $\langle{\bf e},{\bf d}-{\bf e}\rangle\geq 0$, we have
$$\epsilon_{\bf d}^{\rm eff}(\delta)\leq\epsilon^{\rm opt}_{\bf d}(\delta),$$
and this value is sometimes easier to determine. This motivates our strategy for proving existence of uniform expanders:

\begin{lemma}\label{lemmaexuniform} Let $({\bf d}^k)_k$ be an unbounded sequence of dimension vectors. Assume that, for every $0<\delta<1$, the sequence $(\epsilon_{{\bf d}^k}^{\rm eff}(\delta))_k$ can be bounded from below by a positive $\epsilon(\delta)$, that is, $\epsilon_{{\bf d}^k}^{\rm eff}(\delta)\geq\epsilon(\delta)$ for all $k$. Then $({\bf d}^k)_k$ supports a uniform expander.
\end{lemma}

\proof Assume that the condition is fulfilled. Then, for any $k\geq 1$, and for any $\delta\in]0,1[$, we have $\mu({\bf d}^k)-\mu({\bf e})\geq\epsilon(\delta)$ for all $0\not={\bf e}\hookrightarrow{\bf d}^k$ such that $\kappa({\bf e})\leq\delta\cdot\kappa({\bf d}^k)$. By Lemma \ref{lemmaexistence}, there exists a $(\delta,\epsilon(\delta))$-expander representation $V^{(k)}$ of dimension vector ${\bf d}^k$. These clearly form an unbounded family, thus a uniform $(\delta,\epsilon(\delta))$-expander. This proves the lemma. 

\section{The case of generalized Kronecker quivers}\label{sectionkronecker}

Before treating the general case, we use results of \cite{REx} to obtain rather complete results on uniform expansion for generalized Kronecker quivers. So assume $m\geq 3$, and let $Q$ be the quiver with $Q_0=\{1,2\}$ and $m$ arrows from $1$ to $2$. Let ${\bf d}=d_1{\bf 1}+d_2{\bf 2}$ be a dimension vector such that $({\bf d},{\bf d})=d_1^2-md_1d_2+d_2^2<0$. 

\begin{theorem} Under this assumption, the family $(k\cdot{\bf d})_k$ supports a uniform expander with explicit lower bound for $\epsilon_{\bf d}(\delta)$.
\end{theorem}

\proof By \cite[Proposition 3.4]{REx}, we have ${\bf e}\hookrightarrow{\bf d}$ if and only if $\langle{\bf e},{\bf d}-{\bf e}\rangle\geq 0$ if and only if $e_2\geq c_{\bf d}(e_1)$, where the function $c_{\bf d}(x)$ is given by
$$c_{\bf d}(x)=\frac{1}{2}(mx+d_2-\sqrt{(mx-d_2)^2+4x(d_1-x)}).$$
Moreover, the proof of \cite[Theorem 4.3]{REx} shows that, for any $0<x<d_1$, we can write
$$c_{\bf d}(x)=(1+\zeta_{\alpha}(x/d_1))\cdot\frac{d_2}{d_1}\cdot x$$
 for a positive $\zeta_{\alpha}(x/d_1)$, where $\zeta_{\alpha}(x/d_1)$, as a function on $x$, is strictly monotonously decreasing and only depends on $\alpha=d_2/d_1$, and not on ${\bf d}$.\\[1ex]
Let $\mu=\Theta/\kappa$ be a stability for $Q$. We can assume that $\mu({\bf d})=0$, thus, without loss of generality, $\Theta({\bf 1})=d_2$, $\Theta({\bf 2})=-d_1$. Using nonnegativity of ${\bf d}$, ${\bf e}$ and $\kappa$, for any $0\not={\bf e}\hookrightarrow{\bf d}$ we can then estimate
$$-\mu({\bf e})=\frac{-d_2e_1+d_1e_2}{\kappa_1e_1+\kappa_2e_2}\geq\frac{-d_2e_1+d_1c_{\bf d}(e_1)}{\kappa_1e_1+\kappa_2c_{\bf d}(e_1)}=\frac{-d_2e_1+(1+\zeta_\alpha(e_1/d_1))d_2e_1}{\kappa_1e_1+\kappa_2(1+\zeta_\alpha(e_1/d_1))d_2e_1/d_1}=$$
$$=\frac{d_1d_2\cdot\zeta_\alpha(e_1/d_1)}{(\kappa_1d_1+\kappa_2d_2)+\kappa_2d_2\cdot\zeta_{\alpha}(e_1/d_1)}.$$
If $\kappa({\bf e})=\delta\cdot\kappa({\bf d})$, we can estimate
$$\delta\cdot(\kappa_1d_1+\kappa_2d_2)=\kappa_1e_1+\kappa_2e_2>\kappa_1e_1+\kappa_2d_2e_1/d_1=(\kappa_1d_1+\kappa_2d_2)e_1/d_1,$$
and thus $e_1/d_1<\delta$. This allows us to finally bound $-\mu({\bf e})$ from below by a term depending only on $\delta$ (and $\alpha$), namely
$$-\mu({\bf e})\geq \frac{d_1d_2\cdot\zeta_{\alpha}(\delta)}{\kappa({\bf d})+\kappa_2d_2\cdot\zeta_{\alpha}(\delta)}=:\epsilon_{\bf d}(\delta).$$
Application of Lemma \ref{lemmaexuniform} finishes the proof.\\[1ex]
{\bf Remark:} Essentially the same calculation shows that a $(\delta,\epsilon)$-expander in the sense of \cite[Definition 4.1]{REx} is a $(\delta',\epsilon')$-expander representation with respect to the above slope function, for
$$\delta'=\frac{\delta\cdot\kappa_1d_1+\kappa_2d_2}{\kappa_1d_1+\kappa_2d_2},\; \epsilon'=\frac{d_1d_2\cdot\epsilon}{\kappa({\bf d})+\kappa_2d_2\cdot\epsilon}.$$

\section{Construction of uniform expanders for wild quivers}\label{wild}

Before applying Lemma \ref{lemmaexuniform} in the generality of arbitrary wild quivers, we recall properties of the Cartan matrix $C$ and their relation to representations of $Q$, mainly following \cite{Kac,Kac2}.\\[1ex]
We can assume the quiver $Q$ to be connected. The Cartan matrix $C$ of $Q$ is then irreducible (that is, it cannot be transformed to proper block diagonal form by a re-indexing of $Q_0$). It is positive definite if and only if $Q$ is of Dynkin type, positive semi-definite if and only if $Q$ is of Dynkin or extended Dynkin type, and indefinite in all other cases, in which $Q$ is called wild.\\[1ex]
Denote by $n$ the cardinality of $Q_0$. By definition, the matrix $2\cdot E_n-C$ has nonnegative entries, and is still irreducible, thus the Perron-Frobenius theorem applies to it. Together with the symmetry of $C$, this immediately yields:

\begin{lemma}\label{lemfp} Assume that $Q$ is wild and connected. Then $C$ has only real eigenvalues $\lambda_1\leq\ldots\leq\lambda_n$, the smallest one $\lambda_1$ is negative and simple (and thus $\lambda_1<\lambda_2$), and admits an eigenvector ${\bf v}\in\mathbb{R}_{>0}Q_0$.\\
Consequently, the infinum of $\max\{\frac{(C{\bf d})_i}{{\bf d}_i}|i\in Q_0\}$ over ${\bf d}\in\mathbb{R}_{>0}Q_0$ equals $\lambda_1$, and the infinum of $\frac{({\bf x},{\bf x})}{{\bf x}\cdot{\bf x}}$ over ${\bf x}\in {\bf v}^{\perp}$ equals $\lambda_2$.
\end{lemma}

Here we used the standard Euclidean scalar product in $\mathbb{R}Q_0$, for which the coordinate elements ${\bf i}$ form an orthonormal basis, ${\bf i}\cdot{\bf j}=\delta_{ij}$.\\[1ex]
The so-called fundamental domain of $\mathbb{R}_{>0}Q_0$ is the cone of all ${\bf d}\in\mathbb{R}_{>0}Q_0$ for which $({\bf d},{\bf i})\leq 0$ for all $i\in Q_0$ (it is a fundamental domain for the Weyl group action on imaginary roots for the Cartan matrix $C$). By \cite[Lemma 1]{Kac2} and \cite[Theorem 6.1]{Scho}, we have:

\begin{theorem}\label{thmschur} If ${\bf d}$ belongs to the fundamental cone and has connected support, then there exists a Schurian representation of dimension vector ${\bf d}$ (more precisely, there exists a dense subset of $R_{\bf d}(Q)$ consisting of such representations).

For general ${\bf d}$, there exists a Schurian representation of dimension vector ${\bf d}$ if and only if there exists a stable representation of dimension vector ${\bf d}$ for any slope function defined using $\Theta=\{{\bf d},\_\}$.\end{theorem}

We now start with the main construction of a uniform expander. Let $Q$ be wild and connected, and let ${\bf d}\in\mathbb{R}_{>0}Q_0$ be a dimension vector in the interior of the fundamental domain, that is, $({\bf d},{\bf i})<0$ for all $i\in Q_0$.\\[1ex]
Then $\kappa=-({\bf d},\_)$ is a positive functional, we define $\Theta=\{{\bf d},\_\}$, and we denote by $\mu$ the associated slope function. By Theorem \ref{thmschur}, we then know that ${\bf d}$ admits a Schurian representation, thus a $\mu$-stable representation. Note that any positive multiple of ${\bf d}$ gives the same slope function $\mu$, and also admits a $\mu$-stable representation. The family $(k{\bf d})_k$ is our candidate for a uniform expander, and we would like to apply Lemma \ref{lemmaexuniform}. The following lemma provides the key technical estimate.

\begin{lemma}\label{boundbydd} There exists an open neighbourhood $\mathcal{U}\subset\mathbb{R}_{>0}Q_0$ of the ray $\mathbb{R}_{>0}\mathbf{v}$ such that, for all ${\bf d}\in \mathcal{U}$,
there exists a positive constant $C$ with the following property: if ${\bf e}\leq{\bf d}$ fulfills $\langle{\bf e},{\bf d}-{\bf e}\rangle\geq 0$, then
$$\{{\bf d},{\bf e}\}\leq C\cdot\eta\cdot(1-\eta)\cdot({\bf d},{\bf d}),$$
where $\kappa({\bf e})=\eta\cdot\kappa({\bf d})$.
\end{lemma}

\proof We have
$$\{{\bf d},{\bf e}\}=\langle{\bf d},{\bf e}\rangle-\langle{\bf e},{\bf d}\rangle,\; ({\bf d},{\bf e})=
\langle{\bf d},{\bf e}\rangle+\langle{\bf e},{\bf d}\rangle,$$
thus the assumption on ${\bf e}$ can be reformulated as
$$0\leq\langle{\bf e},{\bf d}-{\bf e}\rangle=\langle{\bf e},{\bf d}\rangle-\langle{\bf e},{\bf e}\rangle=\frac{1}{2}({\bf d},{\bf e})-\frac{1}{2}\{{\bf d},{\bf e}\}-\frac{1}{2}({\bf e},{\bf e}),$$
which implies $$\{{\bf d},{\bf e}\}\leq({\bf d},{\bf e})-({\bf e},{\bf e}).$$
By the definition of $\kappa$ and of $\eta$, we can write uniquely $${\bf e}=\eta\cdot{\bf d}+{\bf x}$$
for an element ${\bf x}\in\mathbb{R}Q_0$  such that $$({\bf d},{\bf x})=0,\; -\eta\cdot{\bf d}\leq{\bf x}\leq(1-\eta)\cdot{\bf d}.$$
Then $$({\bf d},{\bf e})=\eta\cdot({\bf d},{\bf d}),\; ({\bf e},{\bf e})=\eta^2\cdot({\bf d},{\bf d})+({\bf x},{\bf x}),$$
and thus
$$\{{\bf d},{\bf e}\}\leq\eta\cdot({\bf d},{\bf d})-\eta^2\cdot({\bf d},{\bf d})-({\bf x},{\bf x})=\eta\cdot(1-\eta)\cdot({\bf d},{\bf d})-({\bf x},{\bf x}).$$\
If $({\bf x},{\bf x})\geq 0$, this estimate proves the claim for $C=1$ without further restrictions on ${\bf d}$. So assume $({\bf x},{\bf x})<0$, and write $({\bf x},{\bf x})=\lambda\cdot{\bf x}\cdot{\bf x}$ for negative $\lambda$.
We derive an estimate for ${\bf x}\cdot{\bf x}$. From the bounds $-\eta d_i\leq x_i\leq(1-\eta)d_i$, we first find $(x_i+\eta d_i)\cdot((1-\eta)d_i-x_i)\geq 0$, and thus
$$x_i^2\leq\eta(1-\eta)d_i^2+(1-2\eta)x_id_i$$
for all $i\in Q_0$. Defining $w_i=-({\bf d},{\bf i})>0$, we then estimate
\begin{eqnarray*}{\bf x}\cdot{\bf x}&=&\sum_{i\in Q_0}\frac{d_i}{w_i}\frac{w_i}{d_i}x_i^2\leq 
\max\{\frac{d_i}{w_i}|i\in Q_0\}\cdot\sum_{i\in Q_0}\frac{w_i}{d_i}x_i^2\leq\\
&\leq&\max\{\frac{d_i}{w_i}|i\in Q_0\}\cdot \sum_{i\in Q_0}\frac{w_i}{d_i}\left(\eta(1-\eta)d_i^2-(1-2\eta)x_id_i\right)=\\
&=&\max\{\frac{d_i}{w_i}|i\in Q_0\}\cdot(\eta(1-\eta)\underbrace{\sum_{i\in Q_0}w_id_i}_{=-({\bf d},{\bf d})}-(1-2\eta)\underbrace{\sum_{i\in Q_0}w_ix_i}_{=-({\bf d},{\bf x})=0})=\\
&=&-\max\{\frac{d_i}{w_i}|i\in Q_0\}\cdot\eta(1-\eta)\cdot({\bf d},{\bf d}).
\end{eqnarray*}
This implies (since $\lambda<0$ by assumption)
$$({\bf x},{\bf x})=\lambda\cdot{\bf x}\cdot{\bf x}\geq -\lambda\cdot\max\{\frac{d_i}{w_i}|i\in Q_0\}\cdot\eta(1-\eta)\cdot({\bf d},{\bf d}),$$
which yields the estimate
$$\{{\bf d},{\bf e})=\eta(1-\eta)\cdot({\bf d},{\bf d})-({\bf x},{\bf x})\leq\underbrace{(1+\lambda\cdot\max\{\frac{d_i}{w_i}|i\in Q_0\})}_{=:C}\cdot\eta(1-\eta)\cdot({\bf d},{\bf d}).$$

Now we use the statements on infima in Lemma \ref{lemfp}. For ${\bf d}$ in a small open neighbourhood $\mathcal{U}$ of $\mathbb{R}_{>0}\mathbf{v}$, the maximum $M:=\max\{\frac{d_i}{w_i}|i\in Q_0\}$ is arbitrarily close to $-\frac{1}{\lambda_1}$. On the other hand, since $({\bf d},{\bf x})=0$, the vector ${\bf x}$ then belongs to an arbitrarily small neighbourhood of the hyperplane $\mathbf{v}^\perp$, thus $\lambda$ is arbitrarily close to $\lambda_2>\lambda_1$. Consequently, the constant $C$ is arbitrarily close to $1-\frac{\lambda_2}{\lambda_1}$, thus in particular positive for small enough $\mathcal{U}$, proving the claim.



\begin{corollary}\label{corex} Under the assumptions of the previous lemma, for all $0<\delta<1$, we have $\epsilon_{\bf d}^{\rm eff}(\delta)\geq C\cdot(1-\delta)>0$. Consequently, the sequence $(k\cdot{\bf d})_k$ supports a uniform expander.
\end{corollary}

\proof Indeed, for $0\not={\bf e}\leq{\bf d}$ such that $\langle{\bf e},{\bf d}-{\bf e}\rangle\geq 0$ and $\kappa({\bf e})\leq\delta\cdot\kappa({\bf d})$, we have (note again that $({\bf d},{\bf d})<0$)
$$\mu({\bf d})-\mu({\bf e})=\frac{\{{\bf d},{\bf e}\}}{({\bf d},{\bf e})}=\frac{\{{\bf d},{\bf e}\}}{\eta\cdot({\bf d},{\bf d})}\geq\frac{C\cdot\eta\cdot(1-\eta)\cdot({\bf d},{\bf d})}{\eta\cdot({\bf d},{\bf d})}=C\cdot(1-\eta)\geq C\cdot(1-\delta),$$
proving the first claim. The second claim follows from Lemma \ref{lemmaexuniform}.

\section{Nonexistence of preprojective/preinjective uniform expanders}\label{nonpp}

To complement our existence result from the previous section, we prove:

\begin{proposition}\label{nopp} There are no uniform expanders consisting entirely of preprojective or preinjective representations.
\end{proposition}

We first recall the definition of preprojective/preinjective representations, in particular limit properties obtained in \cite{DL,PT,RiTame,Ru}. We assume that $Q$ is not of Dynkin type and connected. The quiver $Q$ being acyclic, its path algebra $FQ$ is finite-dimensional, and the regular representation $FQ$ decomposes into indecomposable projective ones $P_i$ for $i\in Q_0$, such that $P_i$ is spanned by all paths starting in $i$. Similarly, the linear dual $(FQ)^*$ decomposes into injective indecomposable representations $I_i$ for $i\in Q_0$. Considering $FQ$ and $(FQ)^*$ as bimodules over $FQ$, we consider the endofunctors $$\tau={\rm Ext}^1_{FQ}(\_,FQ)^*,\; \tau^{-1}={\rm Ext}^1_{FQ}((FQ)^*,\_)$$ on left $FQ$-modules, called (inverse) Auslander-Reiten translation. The $\tau^{-k}P_i$ for $k\geq 0$ and $i\in Q_0$, as well as the $\tau^kI_i$ for $k\geq 0$ and $i\in Q_0$, form infinite series of Schurian representations, called preprojective (resp.~preinjective) indecomposables. One special feature of this family is that every subrepresentation of a preprojective indecomposable $\tau^{-k}P_i$ is isomorphic to a direct sum of preprojective indecomposables of the form $\tau^{-l}P_j$ for $l\leq k$ (and dually for factor representations of preinjective ones).\\[1ex]
The effect of $\tau^{(-1)}$ on dimension vectors is described by the Coxeter transformation $\Phi^{(-1)}\in{\rm GL}(\mathbb{Z}Q_0)$. By \cite{Ri}, the spectral radius $\rho$ of $\Phi$ is at least one (and equal to one if and only if $Q$ is of extended Dynkin type, since $Q$ is assumed not to be of Dynkin type), and it is an eigenvalue of $\Phi$, as is $\rho^{-1}$. For wild quivers $Q$, by \cite{PT}, there exist eigenvectors ${\bf y}^{\pm}\in\mathbb{R}_{>0}Q_0$ of $\Phi$ for the eigenvalues $\rho^{\pm 1}$, respectively, such that
$$\lim_{n\rightarrow\infty}\frac{1}{\rho^n}{\rm\bf dim}(\tau^{-n}P_i)=\lambda^-\cdot{\bf y}^-,\;\lim_{n\rightarrow\infty}\frac{1}{\rho^n}{\rm\bf dim}(\tau^{n}I_i)=\lambda^+\cdot{\bf y}^+$$
for positive reals $\lambda^+,\lambda^-$. For $Q$ of extended Dynkin type, by the definition of the defect function in \cite[Section 1]{DL}, we can still choose ${\bf y}^+={\bf y}^-\in\mathbb{R}_{>0}Q_0$ to be an eigenvector of $\Phi$ for the eigenvalue $1$, such that
$$\lim_{n\rightarrow\infty}\frac{1}{n}{\rm\bf dim}(\tau^{-n}P_i)=\lambda^-\cdot{\bf y}^-,\; \lim_{n\rightarrow\infty}\frac{1}{n}{\rm\bf dim}(\tau^{n}I_i)=\lambda^+\cdot{\bf y}^+$$
for positive reals $\lambda^+,\lambda^-$.\\[1ex]
We can now prove Proposition \ref{nopp}. Assume $(V^{(k)})_{k\geq 1}$ is an unbounded family of dimension vectors of preprojective or preinjective representations. We can choose an unbounded subfamily consisting entirely of one of the types of representations. We only  consider the case of all $V^{(k)}$ being preprojective; the case of preinjective representations can be treated dually. We thus write $$V^{(k)}=\tau^{-n(k)}P_{i(k)},\;\;\; {\bf d}^k={\rm\bf dim}(V^{(k)}).$$
Let $\mu=\Theta/\kappa$ be a slope function for which all $V^{(k)}$ are $\mu$-stable, and fix $\delta\in]0,1[$ and $\epsilon>0$.  Let $i$ be a sink in $Q$. Then $P_i$ is a simple representation, thus any non-zero map from $P_i$ to a preprojective representation $V$ is injective. We have $${\rm Hom}(P_i,V)\simeq V_i,$$ thus $P_i$ embeds into $V$ if ${\rm\bf dim}(V)_i\not=0$. Applying $\tau^{-1}$, we find that $\tau^{-n}P_i$ embeds into a preprojective $V$ if ${\rm\bf dim}(\tau^nV)_i\not=0$. We use the above statements on limits and two properties of the slope function $\mu$: it is continuous, and insensitive to real multiples of a dimension vector. This implies $$\lim_{n\rightarrow\infty}\mu(\tau^{-n}V)=\mu({\bf y}^-)$$ for all preprojective $V$. We can thus choose $n$ large enough so that $$|\mu(\tau^{-n}P_i)-\mu({\bf y}^-)|<\epsilon/2.$$ With the same argument, and again using the growth behaviour of dimension vectors under $\tau^{-1}$, we can also find an index $k$ large enough such that
\begin{itemize} 
\item$(\Phi^{n}{\bf d}^k)_i\not=0$, so that $\tau^{-n}P_i$ embeds into $V^{(k)}$,
\item $\kappa(\tau^{-n}P_i)\leq\delta\cdot\kappa({\bf d}^k)$,
\item $|\mu({\bf d}^k)-\mu({\bf y}^-)|<\epsilon/2$.
\end{itemize}
We then find (since $V^{(k)}$ is $\mu$-stable) that $$\mu(V^{(k)})-\mu(\tau^{-n}P_i)=|\mu(V^{(k)})-\mu(\tau^{-n}P_i)|<\epsilon.$$ This contradicts uniform expansion, proving the proposition. 

\section{Proof of the main result and concluding remarks}\label{concluding}

We can now prove Theorem \ref{main}. If $Q$ is wild, we can choose a dimension vector ${\bf d}$ in an open neighbourhood $\mathcal{U}$ of $\mathbb{R}_{>0}\mathbf{v}$ provided by Lemma \ref{boundbydd}, so that Corollary \ref{corex} provides a uniform expander supported by the sequence $(k\cdot{\bf d})_k$. If $Q$ is Dynkin, there are only finitely many indecomposable, in fact Schurian, representations up to isomorphism, thus there are no unbounded families of representations at all. If $Q$ is extended Dynkin, by \cite[end of Section 3]{DL}, all Schurian representations are preprojective or preinjective, or have dimension vector at most ${\bf y}^-$. Unbounded families thus have to be preprojective/preinjective, and Proposition \ref{nopp} prevents existence of uniform expanders, finishing the proof of the theorem.\\[1ex]
{\bf Remarks:}
\begin{itemize}
\item It is desirable to determine the functions
$$\delta\mapsto\liminf_{k\rightarrow\infty}\epsilon^{\rm eff}_{(k\cdot{\bf d})}(\delta)$$ for a fixed ${\bf d}$ in the interior of the fundamental domain, with respect to the slope function $\mu=-\{{\bf d},\_\}/({\bf d},\_)$ as in Section \ref{wild}. 
\item Our proof of existence of uniform expanders is highly nonconstructive, in that choices of representations outside certain (non-explicit) hypersurfaces have to be made in the proof of Lemma \ref{lemmaexistence}. This also prevents the theorem to generalize to any non-algebraically closed ground field: in choosing expander representations in infinitely many dimension vectors, a priori, coefficients from infinitely many algebraic extensions of the ground field have to be chosen.
\item  In view of our main result, it is tempting to try to define expansion for arbitrary (finite-dimensional) algebras, and to ask whether uniform expansion again characterizes wildness of the algebra. At least this cannot be done in a straightforward way. Namely, the algebras with dense orbit propery exhibited in \cite{CKW} are wild, but there are only finitely many isomorphism classes of Schur representations, thus of stable representations, for any slope stability.
\end{itemize}

\bibliography{expanderrepresentations}{}

\begin{thebibliography}{10}

\bibitem{CKW}
Calin Chindris, Ryan Kinser, and Jerzy Weyman.
\newblock Module varieties and representation type of finite-dimensional
  algebras.
\newblock {\em Int. Math. Res. Not. IMRN}, (3):631--650, 2015.

\bibitem{CB}
William Crawley-Boevey.
\newblock Subrepresentations of general representations of quivers.
\newblock {\em Bull. London Math. Soc.}, 28(4):363--366, 1996.

\bibitem{PT}
J.~A. de~la Pe\~na and M.~Takane.
\newblock Spectral properties of {C}oxeter transformations and applications.
\newblock {\em Arch. Math. (Basel)}, 55(2):120--134, 1990.

\bibitem{DL}
Vlastimil Dlab and Claus~Michael Ringel.
\newblock Indecomposable representations of graphs and algebras.
\newblock {\em Mem. Amer. Math. Soc.}, 6(173):v+57, 1976.

\bibitem{HD}
Lutz Hille and Jos\'e{}~Antonio de~la Pe\~na.
\newblock Stable representations of quivers.
\newblock {\em J. Pure Appl. Algebra}, 172(2-3):205--224, 2002.

\bibitem{Kac}
V.~G. Kac.
\newblock Infinite root systems, representations of graphs and invariant
  theory.
\newblock {\em Invent. Math.}, 56(1):57--92, 1980.

\bibitem{Kac2}
Victor~G. Kac.
\newblock Root systems, representations of quivers and invariant theory.
\newblock In {\em Invariant theory ({M}ontecatini, 1982)}, volume 996 of {\em
  Lecture Notes in Math.}, pages 74--108. Springer, Berlin, 1983.

\bibitem{K}
A.~D. King.
\newblock Moduli of representations of finite-dimensional algebras.
\newblock {\em Quart. J. Math. Oxford Ser. (2)}, 45(180):515--530, 1994.

\bibitem{WW}
Yinan Li, Youming Qiao, Avi Wigderson, Yuval Wigderson, and Chuanqi Zhang.
\newblock On linear-algebraic notions of expansion, 2023.

\bibitem{Lu}
Alexander Lubotzky.
\newblock Expander graphs in pure and applied mathematics.
\newblock {\em Bull. Amer. Math. Soc. (N.S.)}, 49(1):113--162, 2012.

\bibitem{LZ}
Alexander Lubotzky and Efim Zelmanov.
\newblock Dimension expanders.
\newblock {\em J. Algebra}, 319(2):730--738, 2008.

\bibitem{REx}
Markus Reineke.
\newblock Dimension expanders via quiver representations.
\newblock {\em J. Comb. Algebra}, 8(1-2):111--119, 2024.

\bibitem{RiTame}
Claus~Michael Ringel.
\newblock {\em Tame algebras and integral quadratic forms}, volume 1099 of {\em
  Lecture Notes in Mathematics}.
\newblock Springer-Verlag, Berlin, 1984.

\bibitem{Ri}
Claus~Michael Ringel.
\newblock The spectral radius of the {C}oxeter transformations for a
  generalized {C}artan matrix.
\newblock {\em Math. Ann.}, 300(2):331--339, 1994.

\bibitem{Ru}
Alexei Rudakov.
\newblock Stability for an abelian category.
\newblock {\em J. Algebra}, 197(1):231--245, 1997.

\bibitem{Schiffler}
Ralf Schiffler.
\newblock {\em Quiver representations}.
\newblock CMS Books in Mathematics/Ouvrages de Math\'ematiques de la SMC.
  Springer, Cham, 2014.

\bibitem{Scho}
Aidan Schofield.
\newblock General representations of quivers.
\newblock {\em Proc. London Math. Soc. (3)}, 65(1):46--64, 1992.

\end{thebibliography}
\bibliographystyle{plain}


\end{document}